\documentclass[12pt,reqno,a4paper]{amsart}            

\usepackage{hyperref}

\usepackage{scrlayer-scrpage}

\usepackage{scrtime}     

\usepackage{amsmath, amssymb}
\usepackage{amsthm}

\usepackage{accents}     
\usepackage{mathtools}  

\usepackage{bbm}         

\usepackage{mathrsfs}    

\usepackage{graphicx}                

\usepackage[totalwidth=450pt, totalheight=750pt]{geometry}

\newtheorem{theorem}{Theorem}[section]

\newtheorem{proposition}[theorem]{Proposition}
\newtheorem{lemma}[theorem]{Lemma}

\theoremstyle{definition}       
\newtheorem{definition}[theorem]{Definition}

\newtheorem{remark}[theorem]{Remark}

\newtheorem*{remark*}{Remark}

\numberwithin{equation}{section}

\newcommand{\Sec}[1]{Section~\ref{sec:#1}}

\newcommand{\Subsec}[1]{Subsection~\ref{ssec:#1}}

\newcommand{\Thm}[1]{Theorem~\ref{thm:#1}}

\newcommand{\Thmenum}[2]{Theorem~\ref{thm:#1}~(\ref{#2})}

\newcommand{\Lem}[1]{Lemma~\ref{lem:#1}}

\newcommand{\Lemenum}[2]{Lemma~\ref{lem:#1}~(\ref{#2})}

\newcommand{\Prp}[1]{Proposition~\ref{prp:#1}}
\newcommand{\Prps}[2]{Propositions~\ref{prp:#1} and~\ref{prp:#2}}
\newcommand{\Prpss}[3]{Propositions~\ref{prp:#1}, ~\ref{prp:#2} and~\ref{prp:#3}}

\newcommand{\Prpenum}[2]{Proposition~\ref{prp:#1}~(\ref{#2})}

\newcommand{\Rem}[1]{Remark~\ref{rem:#1}}

\newcommand{\Def}[1]{Definition~\ref{def:#1}}

\usepackage{amsmath, amssymb}
\usepackage{amsthm}

\usepackage{accents}     
\usepackage{mathtools}  

\usepackage{bbm}         

\usepackage{mathrsfs}    

\usepackage{xspace}  

\newcommand{\comment}[1]{}   

\usepackage{graphicx}                

\newcommand{\itemref}[1]{\eqref{#1}}

\newcommand{\myparagraph}[1]{\noindent\textbf{#1}}
\newcommand{\abs}[2][{}]{\lvert{#2}\rvert_{{#1}}}    
\newcommand{\abssqr}[2][{}]{\lvert{#2}\rvert^2_{#1}} 
\newcommand{\bigabssqr}[2][{}]{\bigl\lvert{#2}\bigr\rvert^2_{#1}}
\newcommand{\Bigabssqr}[2][{}]{\Bigl\lvert{#2}\Bigr\rvert^2_{#1}}

\newcommand{\normsymb}{\|}
\newcommand{\bignormsymb}[1]{#1\|}

\newcommand{\norm}[2][{}]{\normsymb{#2}\normsymb_{{#1}}}    
\newcommand{\normsqr}[2][{}]{\normsymb{#2}\normsymb^2_{#1}} 
\newcommand{\bignorm}[2][{}]{\bignormsymb{\bigl}{#2}\bignormsymb{\bigr}_{#1}}

\newcommand{\iprod}[3][{}]{\langle{#2},{#3}\rangle_{#1}}  

\newcommand{\set}[2]{\{ \, #1 \, | \, #2 \, \} }      
\newcommand{\bigset}[2]{\bigl\{ \, #1 \, \bigl|\bigr. \, #2 \, \bigr\} }
\newcommand{\Bigset}[2]{\Bigl\{ \, #1 \, \Bigl|\Bigr. \, #2 \, \Bigr\} }

\DeclareMathOperator*{\bigdcup}{\mathaccent\cdot{\bigcup}}
\DeclareMathOperator*{\dcup}   {\mathaccent\cdot\cup}

\newcommand{\map}[3]{ #1 \colon #2 \longrightarrow #3}    

\newcommand{\bd}  {\partial}          

\def\Xint#1{\mathchoice
   {\XXint\displaystyle\textstyle{#1}}%
   {\XXint\textstyle\scriptstyle{#1}}%
   {\XXint\scriptstyle\scriptscriptstyle{#1}}%
   {\XXint\scriptscriptstyle\scriptscriptstyle{#1}}%
   \!\int}
\def\XXint#1#2#3{{\setbox0=\hbox{$#1{#2#3}{\int}$}
     \vcenter{\hbox{$#2#3$}}\kern-.5\wd0}}
\def\XXsum#1#2#3{{\setbox0=\hbox{$#1{#2#3}{\int}$}
     \vcenter{\hbox{$#2#3$}}\kern-.60\wd0}}

\newcommand{\dashint}{\Xint-}   
\newcommand{\avint}{{\textstyle\dashint}}   

\newcommand{\dd}    {\, \mathrm d}    

\DeclareMathOperator{\dom}    {dom}

\DeclareMathOperator{\id}     {id}   

\DeclareMathOperator{\vol}    {vol}

\newcommand{\specsymb} {\sigma} 

\newcommand{\spec}[2][{}]   {\specsymb_{\mathrm{#1}}(#2)}

\newcommand{\eps}{\varepsilon} 

\renewcommand{\phi}{\varphi}   
\renewcommand{\rho}{\varrho}   

\newcommand{\conj}[1]{\overline {#1}}

\newcommand{\R}{\mathbb{R}} 
\newcommand{\C}{\mathbb{C}} 
\newcommand{\N}{\mathbb{N}} 
\newcommand{\Sphere}{\mathbb{S}} 

\newcommand{\1}{\mathbbm 1}                    

\newcommand{\e}{\mathrm e}  

\newcommand{\normder}{\partial_\mathrm{n}}  

\newcommand{\wt}{\widetilde}           

\newcommand{\HS}{\mathscr H}           

\newcommand{\Sobsymb} {\mathsf H} 

\newcommand{\Lsymb}    {\mathsf L}     
\newcommand{\lsymb}    {\ell}          

\newcommand{\Sobspace}[1][1]{\Sobsymb^{#1}}

\newcommand{\Lpspace}[1][p]    {\Lsymb_{#1}}     
\newcommand{\lpspace}[1][p]    {\lsymb_{#1}}     
\newcommand{\Lsqrspace}    {\Lpspace[2]}     
\newcommand{\lsqrspace}    {\lpspace[2]}          

\newcommand{\Lsqr}[2][{}]{\Lsqrspace^{#1}({#2})} 
\newcommand{\lsqr}[2][{}]{\lsqrspace^{#1}({#2})}   

\newcommand{\Sob}[2][1]{\Sobspace [#1]({#2})}         

\newcommand{\err}{\mathrm o}  
\newcommand{\Err}{\mathrm O}

\newcommand{\quadtext}[1]{\quad\text{#1}\quad}
\newcommand{\qquadtext}[1]{\qquad\text{#1}\qquad}

\usepackage{hyperref}
\usepackage{nicefrac}
\usepackage{enumerate}

\newcommand{\dHausdorff}{d_{\mathrm H}}

\newcommand{\HSset}{\mathsf{HS}} 

\newcommand{\BdOpsymb} {\mathcal B}       
\newcommand{\BdOp}[2][{}]{\BdOpsymb_{#1}({#2})}

\newcommand{\CompBdOpsymb}{\mathcal K}
\newcommand{\CompBdOpClasssymb}{\mathsf K}

\newcommand{\OpSet}{\BdOpsymb_{(0,1]}} 
\newcommand{\OpSetHS}[1]{\BdOp[{(0,1]}]{#1}} 
\newcommand{\CompOpSet}{\CompBdOpsymb_{(0,1]}}
\newcommand{\CompOpClass}{\CompBdOpClasssymb_{(0,1]}}
\newcommand{\CompOpSetHS}[1]{\CompOpSet(#1)}

\newcommand{\wtHS}{{\wt \HS}}

\newcommand{\dUNI}{d_{\text{\normalfont uni}}}   
\newcommand{\dQUE}{d_{\text{\normalfont q-uni}}} 
\newcommand{\dSPEC}{d_{\text{\normalfont spec}}}   

\DeclareMathOperator{\Specsymb}{spec}
\renewcommand{\specsymb}{\Specsymb}

\newcommand{\vxeps}{{\eps,v}}
\newcommand{\edeps}{{\eps,e}}

\newcommand{\openint}[2]{(#1,#2)}  
\newcommand{\closedint}[2]{[#1,#2]}  

\newcommand{\VX}{\mathcal G}   
\newcommand{\VXsym}[1][{}]{\VX^{\mathrm{#1}}}  
\newcommand{\VXmax} {\VXsym[max]}   

\begin{document}

\title[Distance between operators acting in different spaces] {A
  distance between operators acting in different Hilbert spaces and
  operator convergence}

\author{Olaf Post}
\address{Fachbereich 4 -- Mathematik,
  Universit\"at Trier,
  54286 Trier, Germany}
\email{olaf.post@uni-trier.de}

\author{Jan Simmer}%
\address{Fachbereich 4 -- Mathematik,
  Universit\"at Trier,
  54286 Trier, Germany}
\email{jan.simmer@yahoo.de}

\ifthenelse{\isundefined \finalVersion} %
{\date{\today, \thistime,  \emph{File:} \texttt{\jobname.tex}}}
{\date{\today}}  



\maketitle



%
\section{Introduction}
\label{sec:intro}
%

The aim of the present article is to give an introduction to the
concept of quasi-unitary equivalence and to define several
(pseudo-)metrics on the space of self-adjoint operators acting
possibly in different Hilbert spaces.  As some of the ``metrics'' do
not fulfil all properties of a metric (e.g.\ some lack the triangle
inequality or the definiteness), we call them ``distances'' here.  To
the best of our knowledge, such distances are treated for the first
time here.  The present article shall serve as a starting point of
further research.

\subsection{Operator convergence in varying Hilbert spaces}
A main motivation for the definition of a distance for operators
acting in different Hilbert spaces is apparent: In many applications,
operators such as a Laplacian $\Delta_\eps \ge 0$ act on a Hilbert space
$\HS_\eps$ that changes with respect to a parameter $\eps$, and one is
interested in some sort of convergence.  Our concept allows to define
a \emph{generalised norm convergence} for the resolvents
$R_\eps=(\Delta_\eps+1)^{-1}$ acting on $\HS_\eps$ towards a resolvent
$R_0=(\Delta_0+1)^{-1}$ acting on $\HS_0$ using identification
operators $\map {J_\eps}{\HS_0}{\HS_\eps}$.  One can first assume that
$J_\eps$ is unitary and that
\begin{subequations}
  \begin{equation}
    \label{eq:que-intro-1}
    \norm{J_\eps R_0 - R_\eps J_\eps} \to 0
  \end{equation}
  as $\eps \to 0$.  In applications (as the one presented in
  \Sec{examples} on shrinking manifolds) it is more convenient to use
  maps $J_\eps$ that are unitary only in an asymptotic sense, i.e.,
  where
  \begin{equation}
    \label{eq:que-intro-2}
    \norm{(\id_{\HS_0}-J_\eps^*J_\eps) R_0} \to 0
    \qquadtext{and}
    \norm{(\id_{\HS_\eps}-J_\eps J_\eps^*) R_\eps} \to 0.
  \end{equation}
\end{subequations}
We call such operators $J_\eps$ \emph{quasi-unitary}, see
\Subsec{que}.  For example, the second estimate
of~\eqref{eq:que-intro-2} means that if $(u_\eps)_\eps$ is a family
with $\norm[\HS_\eps]{(\Delta_\eps+1) u_\eps}=1$, then
$\norm{u_\eps-J_\eps J_\eps^* u_\eps} \to 0$ as $\eps \to 0$.  In our
example, we even have $J^*_\eps J_\eps=\id_{\HS_0}$, and functions in
the range of $\id_{\HS_\eps} - J_\eps J_\eps^*$ do not concentrate at
``negligable'' regions and at high (transversal) modes,
see~\eqref{eq:que-1-ex}--\eqref{eq:que-2-ex}.

We illustrate in \Sec{examples} the abstract theory on generalised
norm resolvent convergence: Consider a family of thin Riemannian
manifolds $X_\eps$ that shrink towards a metric graph $X_0$ (i.e., a
topological graph where the edges are metrically identified with
compact intervals).  We show that the Laplacians on $X_\eps$ converge
in generalised norm resolvent sense to the so-called Kirchhoff
Laplacian on $X_0$.  The example of thin branched manifolds shrinking
towards a metric graph has already been treated in~\cite{post:06} (see
also~\cite{exner-post:05,post:12,exner-post:13} and references
therein).  In this note we use a sightly different proof as we
directly compare the resolvent difference and we do not make use of
the corresponding quadratic forms as in~\cite{post:06,post:12}.  Other
topological perturbations of manifolds such as removing many small
balls are treated in a similar way
in~\cite{khrabustovskyi-post:18,anne-post:pre18}, see also the
references therein.  The concept of generalised norm resolvent
convergence also applies to approximations of metric spaces with a
Laplace-like operator by finite dimensional operators such as graph
approximations of fractals,
see~\cite{post-simmer:18,post-simmer:pre18a,post-simmer:pre18d} for
details.

\subsection{Metrics on sets of operators acting in different Hilbert
  spaces}
When defining a distance between unbounded operators such as
Laplacians, it is convenient to work with the resolvent
$R=(\Delta+1)^{-1}$ where $\Delta$ is an unbounded, self-adjoint and
non-negative operator in a Hilbert space $\HS$.  In particular, we
consider the space of self-adjoint, injective and bounded operators
$R$ with spectrum in $[0,1]$ as space of operators.  In all our
examples, the distance will not change when passing from an operator
$R$ to a unitarily equivalent operator $URU^*$ for a unitary map
$\map U \HS \wtHS$.  The simplest distance we define is
\begin{equation}
  \label{eq:uni.dist}
  \dUNI(R,\wt R)
  := \inf \bigset{\norm{\wt R - UR U^*}}
  {\text{$\map U \HS \wtHS$ unitary}}
\end{equation}
for operators $R$ on $\HS$ and $\wt R$ on $\wtHS$ as above.
If~\eqref{eq:que-intro-1} is fulfilled for some unitary map $J_\eps$,
then $\dUNI(R_\eps,R_0)\to 0$.

From an abstract point of view we could also work with operators in a
\emph{fixed} Hilbert space $\HS$ using an abstract unitary map $U$,
but the identification is in general not natural.  For example, if $R$
and $\wt R$ are both compact, then one can define a unitary map via
$U \psi_k = \wt \psi_k$, where $(\psi_k)_k$ resp.\ $(\wt \psi_k)_k$
are orthonormal bases of eigenfunctions of $R$ resp.\ $\wt R$.  Then
$\norm{\wt R - U R U^*}=\sup_k \abs{\wt \mu_k-\mu_k}$ where $\mu_k$
resp.\ $\wt \mu_k$ denote the corresponding eigenvalues.  This
observation is not very useful in examples, as one needs at least
information on one of the eigenfunction or eigenvalue families.

Later on, we want to use more general maps $\map J \HS \wtHS$ instead
of unitary ones, and allow $J$ to be unitary only ``up to a small
error'', measured e.g.\ by quantities such as $\norm{(\id_\HS-J^*J)R}$
and $\norm{(\id_\wtHS-JJ^*)\wt R}$ (see \Subsec{que} for details).

If $R$ is a compact operator, then more can be said.  Basically, the
different distances defined later on (such as $\dUNI$) depend only on
the spectrum, i.e., the sequence of eigenvalues $(\mu_k)_k$ (ordered
non-increasingly and repeated according to multiplicity).  In particular,
for compact $R$ and $\wt R$, we have $\dUNI(R,\wt R)=0$ if and only if
$R$ and $\wt R$ are unitarily equivalent (this is no longer true for
general $R$ and $\wt R$, see \Rem{not.uni.eq}).

In this article, we only treat \emph{operators}: there is a more
elaborated version of the concept of quasi-unitary equivalence for
closed quadratic forms using not only identification operators
$\map J \HS \wtHS$ on the level of the Hilbert spaces, but also
identification operators $\map {J^1}{\HS^1}{\wtHS^1}$ and
$\map {J^{\prime1}}{\wtHS^1}{\HS^1}$ on the level of the form domains
$\HS^1=\dom \Delta^{1/2}$ and $\wtHS^1=\dom \wt \Delta^{1/2}$.

\subsection{Related works}
There are of course a lot of classical results on operator converence
(and resolvent convergence) for operators acting in a \emph{fixed}
Hilbert space, see e.g.~\cite[Sec.~IV.2]{kato:66}
or~\cite[Ch.~VIII.7]{reed-simon-1}.  The concept of \emph{generalised
  norm resolvent convergence} has already been introduced
by~\cite[Sec.~9.3]{weidmann:00} and is closely related to ours: a
sequence of self-adjoint operators $\Delta_n \ge 0$ converges in
generalised norm resolvent sense to $\Delta_\infty$ if and only if
$R_n=(\Delta_n+1)^{-1}$ and $R_\infty=(\Delta_\infty+1)^{-1}$ are
$\delta_n$-quasi-unitarily equivalent (see \Def{que}) with
$\delta_n \to 0$ as $n \to \infty$, or, equivalently, if
$\dQUE(R_n,R_\infty) \to 0$, see \Def{dist.op-uni}.  B\"ogli recently
proved spectral convergence of spectra and pseudospectra
in~\cite{boegli:17,boegli:18}, we refer also to the extensive list of
references therein.

We compare the different concepts of Kuwae and
Shioya~\cite[Sec.~2]{kuwae-shioya:03} (and related concepts such as
Mosco and $\Gamma$-convergence) generalising \emph{strong} resolvent
convergence, the concept of the generalised norm (and strong)
resolvent convergence in the sense of Weidmann and the results of
B\"ogli in a subsequent paper.  Here, we focus on the definition of
some metrics on a set of operators defined on different Hilbert
spaces, leading to express (operator norm) convergence in metric
terms.  Note that our concept easily allows to define a
\emph{convergence} speed, which in many other works is not treated.

\subsection*{Acknowledgements}
OP would like to thank the organisers of the workshop ``Control theory
of infinite-dimensional systems'' at the Fernuniversit\"at Hagen for
the kind invitation.  We would also like to thank the anonymous
referee for carefully reading our manuscript and helpful comments.

%
\section{
  Distances between operators acting in different Hilbert spaces}
\label{sec:que-definitions}
%

In this section, we introduce a generalisation of norm difference for
operators acting in different Hilbert spaces.

\subsection{A spectral distance for operators acting in different
  Hilbert spaces}
\label{ssec:spec-op-dist}

For a Hilbert space $\HS$, denote by
\begin{equation}
  \label{eq:bd.op.1}
  \OpSetHS \HS
  := \bigset{\map R \HS \HS}
  {R=R^*, \; \ker R=\{0\}, \; \spec R \subset [0,1]}
\end{equation}
the set of all self-adjoint and injective operators with spectrum in
$[0,1]$, i.e., the set of non-negative, self-adjoint and injective
operators with operator norm bounded by $1$.\footnote{We could use any
  other positive number $c>0$ as norm bound, but $1$ makes the
  following estimates simpler and $1$ is also the norm bound in our
  main application where $R=(\Delta+1)^{-1}$ for some non-negative,
  self-adjoint and possibly unbounded operator $\Delta$.}  Moreover,
let $\HSset$ be a set of separable Hilbert spaces of infinite
dimension\footnote{We need to define a fixed set of Hilbert spaces to
  avoid some set-theoretic problems related to self-referencing
  definitions such as ``the set of all sets \dots''. Typically,
  $\HSset$ is s family of Hilbert spaces such as
  $\HSset=\set{\HS_m}{m \in \N \cup \{\infty\}}$ or
  $\HSset=\set{\HS_\eps}{\eps \in [0,1]}$.  Moreover, we assume here
  for simplicity that all Hilbert spaces have infinite dimension.} %
and let
\begin{equation}
  \label{eq:bd.op.set}
  \OpSet
  := \bigcup_{\HS \in \HSset}
  \OpSetHS \HS.
\end{equation}

We first define the following distance function:
\begin{definition}
  \label{def:dist.op-que}
  For $R,\wt R \in \OpSet$ we define the \emph{unitary distance} of
  $R$ and $\wt R$ as in~\eqref{eq:uni.dist}, i.e.,
  \begin{equation*}
    \dUNI(R,\wt R)
    := \inf \bigset{\norm{\wt R - UR U^*}}
      {\text{$U$ unitary}}.
  \end{equation*}
\end{definition}

\begin{proposition}
  \label{prp:dist.op-uni}
  The function $\dUNI$ is a pseudometric on $\OpSet$ (i.e., it
  is a metric execpt for the positive definiteness).  Moreover,
  $\dUNI(R,\wt R)=0$ is equivalent with the fact that there is a
  sequence of unitary operators $\map {U_n}\HS \wtHS$ such that
  \begin{equation}
    \label{eq:quasi-uni-0}
    \norm{\wt R-U_nR U_n^*} \to 0
    \qquad\text{as $n \to \infty$}.
  \end{equation}
  Finally, if $R$ and $\wt R$ are unitarily equivalent, then
  $\dUNI(R,\wt R) = 0$.
\end{proposition}
\begin{proof}
  We have $\dUNI(R,\wt R) \ge 0$, $\dUNI(R,R)=0$,
  $\dUNI(R,\wt R)=\dUNI(\wt R,R)$ and the triangle inequality
  \begin{equation*}
    \dUNI(R_1,R_3)
    \le \dUNI(R_1,R_2) + \dUNI(R_2,R_3).
  \end{equation*}
  follows from
  \begin{multline*}
    \norm{R_3-U_{13}R_1U_{13}^*}
    \le \norm{R_3-U_{23}R_2U_{23}^*} + \norm{U_{23}(R_2-U_{12}R_1U_{12}^*)U_{23}^*}\\
    = \norm{R_3-U_{23}R_2U_{23}^*} + \norm{R_2-U_{12}R_1U_{12}^*}
  \end{multline*}
  using $\map{U_{ij}}{\HS_i}{\HS_j}$ as unitary operators with
  $U_{13}=U_{23}U_{12}$. Taking the infimum over all unitary operators
  $U_{12}$ and $U_{23}$ we obtain the desired inequality.  Note that
  all unitary operators $\map {U_{13}}{\HS_1}{\HS_3}$ can be written as
    $U_{23}U_{12}$, e.g., with $U_{23}=U_{13}U_{12}^*$ for some fixed
    $U_{12}$.  The remaining claims are easily seen.
\end{proof}

\begin{remark}
  \label{rem:not.uni.eq}
  The condition of the two operators $R$ and $\wt R$
  in~\eqref{eq:quasi-uni-0} is closely related to the notion
  \emph{approximate unitary equivalence} defined for the C*-algebras
  generated by $R$ and $\wt R$, cf.~\cite{lin-huaxin:12} and
  references therein.  Note that the unitary orbit of $R$ in
  $\BdOp \wtHS$, i.e., the set $\set{U R U^*}{\text{$\map U \HS \wtHS$
      unitary}}$
  is \emph{not} closed in the operator topology; in particular,
  $\wt R$ and $R$ are not (necessarily) unitarily equivalent if
  $\dUNI(R,\wt R)=0$.  It follows from the next result (see also
  \Prp{dist.op-que}) that such operators must have the same spectrum.
\end{remark}

\begin{proposition}[{\cite[Lem.~A.1]{herbst-nakamura:99}}]
  \label{prp:dist.uni.hausdorff}
  We have
  \begin{equation}
    \label{eq:dist.uni.hausdorff}
    \dHausdorff(\spec R,\spec {\wt R})
    \le \dUNI(R,\wt R),
  \end{equation}
  where $\dHausdorff$ denotes the Hausdorff distance.
\end{proposition}
\begin{proof}
  For a unitary map $\map U \HS \wtHS$, we have
  $\spec {U R U^*}=\spec R$.  In Lemma~A.1
  of~\cite{herbst-nakamura:99} it is proved that
  \begin{equation*}
    \dHausdorff(\spec R,\spec {\wt R})
    \le \norm{U R U^* - \wt R}.
  \end{equation*}
  As $U$ is arbitrary, the assertion follows.
\end{proof}

We restrict now our space of operators to certain compact operators:
Denote by
\begin{equation*}
  \CompOpSetHS \HS
  := \set{R \in \OpSetHS \HS}{\text{$R$ compact}},
\end{equation*}
i.e., the set of compact and injective operators such that
$\spec R \subset [0,1]$.
Moreover, set
\begin{equation}
  \label{eq:bd.comp.set}
  \CompOpSet
  := \bigcup_{\HS \in \HSset}
  \CompOpSetHS \HS.
\end{equation}
For $R \in \CompOpSetHS \HS$, denote by $(\mu_k)_k$ its (discrete
spectrum), ordered in non-increasing order, repeated according to
multiplicity.  Note that $\mu_k \to 0$ as $k \to \infty$.  Denote by
$\Sigma$ the space of all such sequences, i.e.,
\begin{equation}
  \label{eq:def.seq.space}
    \Sigma := \bigset{\mu=(\mu_k)_{k \in \N}}
    {\lim_{k \to \infty} \mu_k=0, \;
      \forall k \in \N \colon 0< \mu_{k+1} \le \mu_k \le 1}
\end{equation}
Denote by $(\psi_k)_k$ a corresponding sequence of orthonormal
eigenfunctions.  As $R$ is injective, $(\psi_k)_k$ is an orthonormal
basis.
Similarly, let $(\wt \mu_k)_k \in \Sigma$ and $(\wt \psi_k)_k$ be the
ordered eigenvalue sequence with corresponding orthonormal basis of
eigenvectors for $\wt R \in \CompOpSetHS \wtHS$.

We set
\begin{equation}
  \label{eq:def.comp.op.class}
  \CompOpClass := \CompOpSet/{\sim}
\end{equation}
where $R \sim \wt R$ if and only if $R$ and $\wt R$ are unitarily
equivalent.

As the class of operators unitarily equivalent with $R$ is actually
determined by the sequence of eigenvalues $(\mu_k)_k$ we have the
following result:
\begin{lemma}
  The map $\map \sigma \CompOpSet \Sigma$ associating to $R$ its
  ordered sequence of eigenvalues $(\mu_k)_k$ descends to a bijective
  map onto the quotient, i.e.,
  $\map {\wt \sigma} \CompOpClass \Sigma$, $[R] \mapsto \sigma(R)$ is
  well-defined and bijective.
\end{lemma}

The main reason why we restrict to the space of compact operators is
that operators with $\dUNI$-distance $0$ are now actually unitarily
equivalent:
\begin{proposition}
  \label{prp:uni.dist}
  We have $\dUNI(R,\wt R)=0$ if and only if $R$ and $\wt R$ are
  unitarily equivalent.  In particular, $\dUNI$ induces a metric on
  $\CompOpClass$.
\end{proposition}
\begin{proof}
  If $\dUNI(R,\wt R)=0$, then $R$ and $\wt R$ are unitarily equivalent
  by~\eqref{eq:dist.que-uni} and \Prpenum{dist.op-que}{dist.op-que.c}.
  If the eigenvalues of $R$ and $\wt R$ are simple, we could also use
  \Prp{dist.uni.hausdorff} to conclude
  $\dHausdorff(\spec R,\spec {\wt R})=0$ and hence
  $\spec R=\spec {\wt R}$.  The simplicity of the spectra implies that
  $R$ and $\wt R$ are actually unitarily equivalent.
\end{proof}

Let us now define a spectral distance respecting also the multiplicity
of the eigenvalues:
\begin{definition}
\label{def:spec.dist}
  For $R, \wt R \in \CompOpSet$ denote by
  \begin{equation*}
    \dSPEC(R, \wt R)
    := \sup_{k \in \N} \abs{\mu_k-\wt\mu_k}
  \end{equation*}
  the \emph{(multiplicity respecting) distance of the spectra}.
\end{definition}

We have some simple consequences:
\begin{lemma}
\label{lem:spec.dist}
  \indent
  \begin{enumerate}[(a)]
  \item
    \label{spec.dist.a}
    The supremum in \Def{spec.dist} is actually a maximum.
  \item
    \label{spec.dist.c}
    $\dSPEC$ defines a pseudometric on $\CompOpSet$.
  \item
    \label{spec.dist.d}
    We have $\dSPEC(R, \wt R)=0$ if and only if $R$ and $\wt R$ are
    unitarily equivalent.  In particular, $\dSPEC$ induces a metric on
    $\CompOpClass$.
  \end{enumerate}
\end{lemma}
\begin{proof}
  \itemref{spec.dist.a}~This is clear as the sequences are monotonely
  decreasing and converge to $0$.  \itemref{spec.dist.c}~As the right
  hand side in the definition of $\dSPEC$ is the supremum norm, the
  claim is standard. \itemref{spec.dist.d}~If
  $\dSPEC(R,\wt R)=0$ then $\mu_k=\wt \mu_k$ for all
  indices $k$.  Define a unitary map by $U\psi_k=\wt \psi_k$, then
  $\wt R=URU^*$.
\end{proof}

\begin{proposition}
  \label{prp:op-dist-spec}
  For $R, \wt R \in \CompOpSet$, we have
  \begin{equation*}
    \dSPEC(R,\wt R)
    \ge \dUNI(R,\wt R).
  \end{equation*}
\end{proposition}
\begin{proof}
  Let $\map U \HS \wtHS$ be the unitary map given by
  $U\psi_k=\wt \psi_k$, then it is easily seen that
  \begin{equation*}
    \dSPEC(R,\wt R)
    =\abs{\wt \mu_{k_0}-\mu_{k_0}}
    =\norm{(\wt R - URU^*)\wt \psi_{k_0}}
    =\norm{\wt R - URU^*},
  \end{equation*}
  where the maximum is achieved at $k_0$.  As $U$ is unitary, the
  inequality follows by the definition of $\dUNI(R,\wt R)$ via an
  infimum over all unitary maps $\map U \HS \wtHS$.
\end{proof}

Passing to the sequence space $\Sigma$, we define
\begin{equation*}
  \dHausdorff(\mu,\wt \mu)
  := \dHausdorff(\set{\mu_k}{k \in \N}, \set{\wt\mu_k}{k \in \N}).
\end{equation*}
This is actually only a pseudometric as the multiple appearance of a value
in a \emph{sequence} $\mu=(\mu_k)_k$ is not detected in the \emph{set}
$\set{\mu_k}{k \in \N}$.  Using the symbols $\dUNI$ and $\dSPEC$ also
for the induced metrics on $\Sigma$ (see~\eqref{eq:def.seq.space}), we
have
\begin{equation}
  \label{eq:three.metrics}
  \dHausdorff(\mu,\wt \mu)
  \le \dUNI(\mu, \wt \mu)
  \le \dSPEC(\mu,\wt \mu)
\end{equation}
combining \Prps{dist.uni.hausdorff}{op-dist-spec}.

Note that the metric space $(\Sigma,\dSPEC)$ is not complete, choose
e.g.\ the sequence $(\mu^{(n)})_n$ with
$\mu^{(n)}=(1/(kn))_{k \in \N}$, then
$\dSPEC(\mu^{(n)},\mu^{(m)})=\abs{1/n-1/m} \to 0$ as $m,n \to \infty$,
i.e., $(\mu^{(n)})_n$ is a Cauchy sequence but the limit $0=(0)_k$ is
not in $\Sigma$.  Similarly, $(\Sigma,\dHausdorff)$ is not complete
(as the set $\set{\mu_k}{k \in \N}$ is not closed in $\R$).  It can be
seen similarly that $(\Sigma,\dUNI)$ is not complete.

\subsection{Quasi-unitary equivalence}
\label{ssec:que}

We now want to weaken the condition that $U$ is unitary in
\Def{dist.op-que} and use a slightly more general concept.  We define
the correspondent distance in \Subsec{op-dist}.
\begin{definition}
  \label{def:que}
  \begin{subequations}
    \label{eq:quasi-uni}
    Let $\delta \ge 0$.  Moreover, let $R \in \OpSetHS \HS$ and
    $\wt R \in \OpSetHS \wtHS$.
      We say that, $R$ and $\wt R$ are
      \emph{$\delta$-quasi-unitarily equivalent}, if there are bounded
      operators $\map J \HS {\wt \HS}$ and $\map {J'}{\wt \HS}\HS$
      such that
      \begin{align}
        \label{eq:quasi-uni.a}
        &\norm J \le 1+\delta,&
                                \norm {J'} \le 1+\delta,&\\
        \label{eq:quasi-uni.b}
        &\norm{J'-J^*} \le \delta,&\\
        \label{eq:quasi-uni.c}
        &\norm{(\id_\HS - J'J)R} \le\delta,&
                                             \norm{(\id_\wtHS - JJ')\wt R}
                                             \le \delta,&\\
        \label{eq:quasi-uni.d}
        &\norm{JR - \wt R J} \le\delta,&
                                         \norm{J'\wt R - R J'}
                                         \le \delta.&
      \end{align}
      We call $J$ and $J'$ \emph{identification operators} and
      $\delta$ the \emph{error}.

  \end{subequations}
\end{definition}
Actually, some conditions follow from others with possibly different
$\delta$, see e.g.~the next lemma; we have included all of them in the
above definition to make them \emph{symmetric} with $R$ and $\wt R$.

Obviously, if $\delta=0$
in~\eqref{eq:quasi-uni.a}--\eqref{eq:quasi-uni.c} then $J$ is unitary,
and~\eqref{eq:quasi-uni.d} is equivalent to the norm estimate
$\norm{\wt R-J^*RJ} \le \delta$.  In particular, $0$-quasi-unitary
equivalence is just unitary equivalence.

For example, we have the following simple facts:
\begin{lemma}
  \label{lem:simplify-que}
  \indent
  \begin{enumerate}[(a)]
  \item Assume that $\norm J \le 1+\delta$ and~\eqref{eq:quasi-uni.b}
    hold, then $\norm {J'} \le 1+2\delta$.
  \item Assume that $\norm {JR-\wt RJ}\le \delta$ and~\eqref{eq:quasi-uni.b}
    hold, then $\norm {J'\wt R-RJ'} \le 3\delta$.
  \item If $J'=J^*$ then $\norm{J}=\norm{J'}$ and
    $\norm{JR-\wt RJ}=\norm{J'\wt R-RJ'}$, i.e., only one of the
    estimates in~\eqref{eq:quasi-uni.a} and~\eqref{eq:quasi-uni.d} is
    enough to ensure $\delta$-quasi-unitary equivalence.
  \item If $J$ is unitary, then $R$ and $\wt R$ are $\delta$-quasi
    unitarily equivalent (with unitary $J$) if and only if
    $\norm{\wt R - JRJ^*}\le \delta$.

  \item
    \label{simplify-que.e}
    If $R_1$ and $R_2$ are $\delta_{12}$-quasi-unitarily equivalent
    and $R_2$ and $R_3$ are $\delta_{23}$-quasi-unitarily equivalent,
    then $R_1$ and $R_3$ are $\delta_{13}$-quasi-unitarily equivalent
    with
    $\delta_{13}=\Phi(\delta_{12},\delta_{23})=\Err(\delta_{12})+\Err(\delta_{23})$,
    where $\Phi$ is defined in~\eqref{eq:def.phi}.

  \item
    \label{simplify-que.g}
    If $R$ and $\wt R$ are $\delta$-quasi-unitarily equivalent with
    $\delta \in [0,1]$ and with identification operators $J$ and $J'$
    then $R$ and $\wt R$ are $3\delta$-quasi-unitarily equivalent with
    identification operators $J$ and $J^*$.  In particular, we can
    assume without loss of generality that $J'=J^*$ in \Def{que}.
  \end{enumerate}
\end{lemma}
\begin{proof}
  The first four assertions are obvious.
  \itemref{simplify-que.e}~The transitivity of quasi-unitary
  equivalence can be seen as follows (see
  also~\cite[Thm.~4.2.5]{post:12} and~\cite{post-simmer:pre18d}):
  Denote by $\map {J_{12}}{\HS_1}{\HS_2}$ and
  $\map {J_{21}}{\HS_2}{\HS_1}$ the identification operators for
  $R_1 \in \OpSetHS {\HS_1}$ and $R_2 \in \OpSetHS {\HS_2}$, and
  similarly, denote by $\map {J_{23}}{\HS_2}{\HS_3}$ and
  $\map {J_{32}}{\HS_3}{\HS_2}$ the identification operators for
  $R_2\in \OpSetHS {\HS_2}$ and $R_3 \in \OpSetHS {\HS_3}$.

  We define the identification operators for $R_1$ and $R_3$ by
  $J_{13}:= J_{23}J_{12}$ and $J_{31}:=J_{21}J_{32}$. Then
  \begin{equation*}
    \norm{J_{13}}
    =\norm{J_{23}J_{12}}
    \le(1+\delta_{23})(1+\delta_{12})
    =1+(\delta_{12}+\delta_{23}) + \delta_{12}\delta_{23}
  \end{equation*}
  and similarly for $\norm{J_{31}}$.
  Inequality~\eqref{eq:quasi-uni.b} follows from
  \begin{multline*}
    \norm{J_{13}^*-J_{31}}
    \le \norm{J_{12}^*(J_{23}^*-J_{32})} + \norm{(J_{12}^*-J_{21})J_{32}}\\
    \leq (1+\delta_{12})\delta_{23} + \delta_{12}(1+\delta_{23})
    = \delta_{12}+\delta_{23} + 2 \delta_{12}\delta_{23}.
  \end{multline*}
  The first inequality in~\eqref{eq:quasi-uni.c} is also satisfied
  because
  \begin{multline*}
    \norm{(\id_{\HS_1} - J_{31} J_{13})R_1}
    \le \norm{(\id_{\HS_1} - J_{21} J_{12})R_1}
      + \norm{J_{21}(J_{12}R_1-R_2J_{12})}\\
      + \norm{J_{21}(\id_{\HS_2}-J_{32}J_{23}) R_2 J_{12}}
      + \norm{J_{21}J_{32}J_{23}(R_2 J_{12}-J_{12}R_1)}\\
    \le \delta_{12}+ (1+\delta_{12})\delta_{12}
      + (1+\delta_{12})^2\delta_{23}
      + (1+\delta_{12})(1+\delta_{23})^2\delta_{12}\\
    = 3\delta_{12}+ \delta_{23}
    + 4\delta_{12}\delta_{23}
    + \delta_{12}^2
    + \delta_{12}^3
      + 3\delta_{12}^2\delta_{23}
      + \delta_{12}\delta_{23}^2
      + \delta_{12}^2 \delta_{23}^2
  \end{multline*}
  and similarly we have
  \begin{multline*}
    \norm{(\id_{\HS_3} - J_{13} J_{31})R_3}
    \le 3\delta_{23}+ \delta_{12}
    + 4\delta_{12}\delta_{23}
    + \delta_{23}^2
      + 3\delta_{12} \delta_{23}^2
      + \delta_{12}^2 \delta_{23}
    + \delta_{23}^3
      + \delta_{12}^2\delta_{23}^2.
  \end{multline*}
  For the first inequality of~\eqref{eq:quasi-uni.d} we estimate
  \begin{multline*}
    \norm{J_{13}R_1 - R_3J_{13}}
    \le \norm {J_{23}(J_{12}R_1-R_2J_{12})} + \norm{(J_{23}R_2-R_3J_{23})J_{12}}\\
    \le (1+\delta_{23})\delta_{12}+\delta_{23}(1+\delta_{12})
    =\delta_{12}+\delta_{23} + 2\delta_{12}\delta_{23}
  \end{multline*}
  and similarly for the second inequality of~\eqref{eq:quasi-uni.d}.
  In particular, the desired estimate holds if we define
  $\delta_{13}:=\Phi(\delta_{12},\delta_{23})$ with
  \begin{equation}
    \label{eq:def.phi}
    \Phi(a,b)
    := 3(a+b) + (a+b)^2 + 2ab + (a+b)^3 + a^2b^2.
  \end{equation}
  Note that $\Phi(a,b)\le 12(a+b)$ if $a,b \in [0,1]$.

  \itemref{simplify-que.g}~We have
  \begin{equation*}
    \norm{(\id_\HS - J^*J)R}
    \le \norm{(\id_\HS - J'J)R} + \norm{J'-J^*}\norm{J}\norm{R}
    \le \delta(2+\delta)
    \le 3\delta
  \end{equation*}
  provided $\delta \in [0,1]$, and similarly for the second inequality
  of~\eqref{eq:quasi-uni.c}.  For the second inequality
  of~\eqref{eq:quasi-uni.d} we have
  \begin{equation*}
    \norm{J^*\wt R-RJ^*}
    \le \norm{(J^*-J')R} + \norm{J'\wt R - RJ'} + \norm{R(J'-J^*)}
    \le 3\delta.   \qedhere
  \end{equation*}
\end{proof}

\subsection{Consequences of quasi-unitary equivalence}
\label{ssec:consequences-que}

Let us cite here some consequences of quasi-unitary equivalence; for
details we refer to~\cite[Ch.~4]{post:12}
and~\cite{post-simmer:pre18d}.  Note that there, we applied
quasi-unitary equivalence to the resolvents $R=(\Delta+1)^{-1}$ and
$\wt R=(\wt\Delta+1)^{-1}$ of two non-negative and self-adjoint
operators $\Delta$ and $\wt\Delta$.
\begin{theorem}
  \label{thm:que}
  \indent
  \begin{enumerate}[(a)]
  \item
    \label{que.a}
    \myparagraph{Convergence of operator functions:} Let
    $\map \phi {[0,1]} \C$ be a continuous function then there are
    functions $\eta_\phi(\delta),\eta_\phi'(\delta) \to 0$ as
    $\delta \to 0$ depending only on $\phi$ such that
    \begin{equation*}
      \norm{\phi(\wt R)J - J \phi(R)}
      \le \eta_\phi (\delta)
      \quadtext{and}
      \norm{\phi(\wt R) - J \phi(R)J'}
      \le \eta'_\phi(\delta)
    \end{equation*}
    for all operators $R$ and $\wt R$ being $\delta$-quasi-unitarily
    equivalent (with identification operators $J$ and $J'$).

    If $\phi$ is holomorphic in a neighbourhood of $[0,1]$ then we can
    choose $\eta_\phi(\delta)=C_\phi \delta$ and similarly for
    $\eta'_\phi$.
  \item \myparagraph{Convergence of spectra:} Let
    $R \in \OpSetHS \HS$, then there is a function $\eta$ with
    $\eta(\delta) \to 0$ as $\delta \to 0$ depending only on $R$ such
    that
    \begin{equation*}
      \dHausdorff(\spec R,\spec{\wt R})
      \le \eta(\delta)
    \end{equation*}
    for all $\wt R$ being $\delta$-quasi-unitarily equivalent with
    $R$.  Here, $\dHausdorff$ denotes the Hausdorff distance of the
    two spectra.  A similar assertion holds for the \emph{essential}
    spectra.

  \item \myparagraph{Convergence of discrete spectra:} Let $\mu$
    be an eigenvalue of $R$ with multiplicity $m$, then there is
    $\delta_0>0$ and there exist $m$ eigenvalues $\wt \mu_j$
    of $\wt R$ ($j=1,\dots, m$, not all necessarily distinct) such
    that $\abs{\mu-\wt\mu_j} \le C \delta$ if
    $\delta \in [0,\delta_0]$, where $C$ is a universal constant
    depending only on $\mu$ and its distance from the remaining
    spectrum of $R$.

  \item \myparagraph{Convergence of eigenfunctions:} Let $\mu$ be
    a simple\footnote{This assumption is for simplicity only.}
    eigenvalue with eigenfunction $\psi$, and denote by $\wt \psi$ the
    eigenfunction associated with $\wt \mu$, then
    $\norm[\wtHS]{\wt \psi - J\psi} \le C' \delta$.  Here, $C'$ is
    again universal constants depending only on $\mu$ and its
    distance from the remaining spectrum of $R$.
  \end{enumerate}
\end{theorem}
The proof of the assertions and more details can be found
in~\cite[Ch.~4]{post:12}.
\begin{remark}
  We say that $\Delta$ and $\wt\Delta$ are
  \emph{$\delta$-quasi-unitarily equivalent} if $R$ and $\wt R$ are.
  \begin{enumerate}[(a)]
  \item
    As an example in~\Thmenum{que}{que.a} one can choose
    $\phi_t((\lambda+1)^{-1})=\e^{-t\lambda}$ and one obtains operator
    estimates for the \emph{heat} or \emph{evolution operator} of
    $\Delta$ and $\wt \Delta$.  In this case, one can use the
    \emph{holomorphic functional calculus} and give a precise estimate
    on $\eta'_{\phi_t}(\delta)$, (namely
    $\eta'_{\phi_t}(\delta)= (16/t+5)\delta$, and we have
    \begin{equation*}
      \bignorm{\e^{-t \wt \Delta} - J' \e^{-t \Delta} J}
      \le \Bigl(\frac {16}t + 5\Bigr)\delta
    \end{equation*}
    provided $\Delta$ and $\wt \Delta$ are $\delta$-quasi-unitarily
    equivalent (or equivalently, $\dQUE(R,\wt R) < \delta$),
    see~\cite{post-simmer:pre18d} for details.  Such an estimate might
    be of interest in control theory.

  \item It suffices in~\itemref{que.a} of this remark that $\phi$ is
    only continuous on $[0,1] \setminus \spec R$.  One can then show
    norm estimates also for \emph{spectral projections}.
  \end{enumerate}
\end{remark}

\subsection{A distance arising from quasi-unitary equivalence}
\label{ssec:op-dist}

We now use the concept of quasi-unitary equivalence to define another
distance function:
\begin{definition}
  \label{def:dist.op-uni}
  For $R,\wt R \in \OpSet$ we define the \emph{quasi-unitary distance}
  of $R$ and $\wt R$ by
  \begin{equation*}
    \dQUE(R,\wt R)
    := \inf \bigset{\delta \ge 0}
    {\text{$R$ and $\wt R$ are $\delta$-quasi-unitarily equivalent}}.
  \end{equation*}
\end{definition}
Clearly, we have
\begin{multline*}
  \dQUE(R,\wt R)
  = \inf \Bigset{
    \max\bigl\{\norm J-1, \norm{J'}-1, \norm{J^*-J'},
    \norm{(\id_\HS - J'J)R},  \norm{(\id_\wtHS - JJ')\wt R},\\
    \norm{JR - \wt R J}, \norm{J'\wt R - R J'}\bigr\}}
  {\text{$\map J \HS \wtHS$, $\map{J'}\wtHS \HS$ bounded}}.
\end{multline*}
Obviously, we have (using also \Prp{op-dist-spec})
\begin{equation}
  \label{eq:dist.que-uni}
  \dQUE(R,\wt R)
  \le \dUNI(R,\wt R)
  \le \dSPEC(R,\wt R)
\end{equation}
as in the definition of $\dUNI$, we only use unitary maps $J$ instead
of general ones in the definition of $\dQUE$.

The function $\dQUE$ has the following properties:
\begin{proposition}
  \label{prp:dist.op-que}
  Let $R,\wt R \in \OpSet$.
  \begin{enumerate}[(a)]
  \item
    \label{dist.op-que.a}
    We have $\dQUE(R,\wt R) \ge 0$, $\dQUE(R,R)=0$,
    $\dQUE(R,\wt R)=\dQUE(\wt R,R)$ and
    \begin{equation*}
      \dQUE(R_1,R_3)
      \le \Phi\bigl(\dQUE(R_1,R_2),\dQUE(R_2,R_3)\bigr)
    \end{equation*}
    where $\Phi(a,b)=3(a+b)+\err(a)+\err(b)$ is defined
    in~\eqref{eq:def.phi}.
  \item
    \label{dist.op-que.b}
    If $\dQUE(R,\wt R)=0$ then the essential spectra of $R$ and
    $\wt R$ agree.  Also the discrete spectra agree and have the same
    multiplicity.
  \item
    \label{dist.op-que.c}
    If $R$ and $\wt R$ are compact, i.e., $R, \wt R \in \CompOpSet$,
    then $\dQUE(R,\wt R)=0$ if and only if $R$ and $\wt R$ are
    unitarily equivalent.  In particular, $\dQUE$ induces a metric on
    $\CompOpClass$ and hence on $\Sigma$
    (see~\eqref{eq:def.seq.space}).
  \end{enumerate}
\end{proposition}
\begin{remark*}
  The function $\dQUE$ in~\itemref{dist.op-que.a} is sometimes
  referred to as \emph{semi(pseudo)metric} with relaxed triangle
  equation, and $\Phi$ is referred to as \emph{triangle function}.  As
  $\Phi(a,b) \le 12(a+b)$ if $a,b \in [0,1]$, and since
  $\dQUE(R,\wt R)\le \dSPEC(R,\wt R)\le 1$, $\dQUE$ is a
  semi(pseudo)metric with $12$-relaxed triangle inequality.
\end{remark*}
\begin{proof}
  \itemref{dist.op-que.a}~For the first assertion, only the relaxed
  triangle inequality is non-trivial; but it follows easily from the
  transitivity in \Lemenum{simplify-que}{simplify-que.e}.

  \itemref{dist.op-que.b}~For the second assertion, note that there
  are bounded operators $\map{J_n}\HS\wtHS$ and $\map{J_n'}\wtHS \HS$
  such that~\eqref{eq:quasi-uni.a}--\eqref{eq:quasi-uni.d} are
  fulfilled for some sequence $\delta_n \to 0$.  In particuluar, $R$
  and $\wt R$ are $\delta_n$-quasi-unitarily equivalent.  It follows
  from \Thm{que} that
  $\dHausdorff(\spec R,\spec {\wt R}) \le \eta(\delta_n)$ for any $n$,
  hence the Hausdorff distance is $0$.  The same is true for the
  essential spectra, and also for the discrete spectrum (including
  multiplicity).  In particular, operators with $\dQUE(R,\wt R)=0$ have
  the same essential and discrete spectrum (the latter even with the
  same multiplicity).

  \itemref{dist.op-que.c}~If $R$ and $\wt R$ are compact, then the
  essential spectrum of both is $\{0\}$ (which is not an eigenvalue)
  and there are orthonormal bases $(\psi_n)_n$ and $(\wt \psi_n)_n$ of
  eigenfunctions of $R$ and $\wt R$, respectively.  Then
  $U\psi_n=\wt \psi_n$ defines a unitary map from $\HS$ to $\wtHS$
  such that $\wt R=URU^*$.
\end{proof}

Due to \Prpenum{dist.op-que}{dist.op-que.c} we can again define a
(semi)metric on $\Sigma$ via $\dQUE(\mu,\wt\mu)=\dQUE(R,\wt R)$ if
$\mu=(\mu_k)_k$ is the sequence of eigenvalues of $R$ and similarly
for $\wt \mu$ and $\wt R$.  Moreover, we have
\begin{equation}
  \label{eq:three.metrics'}
  \dQUE(\mu,\wt \mu)
  \le \dUNI(\mu, \wt \mu)
  \le \dSPEC(\mu,\wt \mu)
\end{equation}
(see~\eqref{eq:three.metrics}).  We will investigate the structure of
$\Sigma$ with respect to these metrics and related questions in a
forthcoming publication.  It is in particular of interest to express
$\dUNI(\mu,\wt \mu)$ and $\dQUE(\mu,\wt \mu)$ directly in terms of the
sequences $\mu$ and $\wt \mu$.

%
\section{Laplacians on thin branched manifolds shrinking towards a
  metric graph}
\label{sec:examples}
%

Let us present our main example here, the Laplacian on a manifold that
shrinks to a metric graph.  Our result holds also for non-compact
Riemannian manifolds and metric graphs under some uniform conditions.
We would like to stress that this example has already been treated
in~\cite{post:06} (see also~\cite{exner-post:05,post:12,exner-post:13}
and references therein).  However, here we present a sightly different
proof as we directly compare the resolvent difference and we do not
make use of the corresponding quadratic forms as
in~\cite{post:06,post:12}.

\subsection{Metric graphs}

Let $(V,E,\bd)$ be a discrete oriented graph, i.e., $V$ and $E$ are at
most countable sets, and $\map \bd E {V \times V}$,
$e \mapsto (\bd_-e,\bd_+e)$ is a map that associates to an edge its
\emph{initial} ($\bd_-e$) and \emph{terminal} ($\bd_+ e)$ vertex.  We
set $E_v^\pm := \set{e \in E}{\bd_\pm e=v}$ and
$E_v := E_v^+ \dcup E_v^-$ (the set of incoming/outgoing resp.\
adjacent edges at $v$).  A \emph{metric} graph is given by $(V,E,\bd)$
together with a map $\map \ell E {\openint 0 \infty}$,
$e \mapsto \ell_e>0$, where we interprete $\ell_e$ as the
\emph{length} of the edge $e$.  In particular, we set
\begin{equation*}
  M_e := \closedint 0 {\ell_e} \quadtext{and}
  M := \bigdcup_{e \in E} M_e / \Psi,
\end{equation*}
where $\Psi$ identifies the end points of the intervals $M_e$
according to the graph structure, i.e.,
\begin{equation*}
  \map \Psi {\bigdcup_{e \in E} \bd M_e} V,
  \qquad
  \begin{cases}
    0 \in \bd M_e \mapsto \bd_-e \in V\\
    \ell_e \in \bd M_e \mapsto \bd_+e \in V
  \end{cases}
\end{equation*}

Any point in $M$ not being a vertex after the identification is
uniquely determined by $e \in E$ and $s_e \in M_e$.  In the sequel, we
often omit the subscript and write $s \in M_e$.

To avoid some technical complications, we assume that
\begin{equation}
  \label{eq:len.bdd}
  \ell_0 := \inf_{e \in E} \ell_e > 0.
\end{equation}
Then $M$ becomes a metric space by defining the distance $d(x,y)$ of
two points $x,y \in M$ as the length of the shortest path in $M$ (the
path may not be unique).  We also have a natural measure on $M$
denoted by $\dd s$, given by the sum of the Lebesgue measures
$\dd s_e$ on $M_e$ (up to the boundary points, a null set).

As Laplacian on $M$, we define $\Delta_M f$ via
$(\Delta_M f)_e=-f_e''$ with $f$ in
\begin{equation}
  \label{eq:dom.lap.mg}
  \dom \Delta_M
  = \bigset{f \in \bigoplus_{e \in E} \Sob[2]{M_e}}
  {\text{$f$ continuous, } \forall v \in V \colon \sum_{e \in E_v} f_e'(v)=0}.
\end{equation}
Here, $f_e'(v)=f_e'(\bd_+ e)$ if $v=\bd_+e$ and $f_e'(v)=-f_e'(0)$ if
$v=\bd_-e$ denotes the derivative of $f$ along the edge $e$ towards
the vertex.  It can be shown that this operator is self-adjoint,
see~\cite{post:12} and references therein for details.  Because of the
sum condition on the derivatives, this operator is also called
\emph{Kirchhoff Laplacian} on $M$.  We later on write $X_0=M$ and
$\Delta_0$ for the Kirchhoff Laplacian $\Delta_M$.

\subsection{Thin branched manifolds (``fat graphs'')}
\label{ssec:graph-likel-mfd}

Let $X_0=M$ be a metric graph.  Let us now describe a family of
manifolds $X_\eps$ shrinking to $X_0$ as $\eps \to 0$.  We will show
that a suitable Laplacian on $X_\eps$ (actually, the Neumann
Laplacian, if $\bd X_\eps \ne \emptyset$) converges to the Kirchhoff
Laplacian $\Delta_0$ on $X_0$

According to the metric graph $X_0$ we associate a family
$(X_\eps)_\eps$ of smooth Riemannian manifolds of dimension
$m+1 \ge 2$ for small $\eps>0$.  We call $X_\eps$ a \emph{thin
  branched manifold}, if the following holds:
\begin{itemize}
\item We have a decomposition
  \begin{equation}
    \label{eq:decomp}
    X_\eps
    = \bigcup_{e \in E} X_\edeps \cup \bigcup_{v \in V} X_\vxeps,
  \end{equation}
  where $X_\edeps$ and $X_\vxeps$ are compact Riemannian manifolds
  with boundary, $(X_\edeps)_{e \in E}$ and $(X_\vxeps)_{v \in V}$ are
  disjoint and
  \begin{equation*}
    X_\edeps \cap X_\vxeps
    =
    \begin{cases}
      \emptyset,& e \notin E_v\\
      Y_\edeps, & e \in E_v,
    \end{cases}
  \end{equation*}
  where $Y_\edeps$ is a compact Riemannian manifold of dimension $m$,
  isometric to a Riemannian manifold $(Y_e, \eps^2 h_e)$.
\item The \emph{edge neighbourhood} $X_\edeps$ is isometric to a
  cylinder
  \begin{equation*}
    X_\edeps \cong M_e \times Y_\edeps,
  \end{equation*}
  i.e., $X_\edeps = M_e \times Y_e$ as manifold with metric $g_\edeps
  = \dd s^2 + \eps^2 h_e$.  Recall that $M_e=\closedint 0 {\ell_e}$.
\item The \emph{vertex neighbourhood} $X_\vxeps$ is isometric to
  \begin{equation*}
    X_\vxeps \cong \eps X_v,
  \end{equation*}
  i.e., $X_\vxeps = X_v$ as manifold with metric $g_\vxeps = \eps^2
  g_v$, where $(X_e,g_v)$ is a compact Riemannian manifold.  In other
  words, $X_\vxeps$ is $\eps$-homothetic with a fixed Riemannian
  manifold $(X_v,g_v)$.
\end{itemize}
\begin{figure}[h]
  \centering
  \begin{picture}(0,0)%
    \includegraphics[scale=0.5]{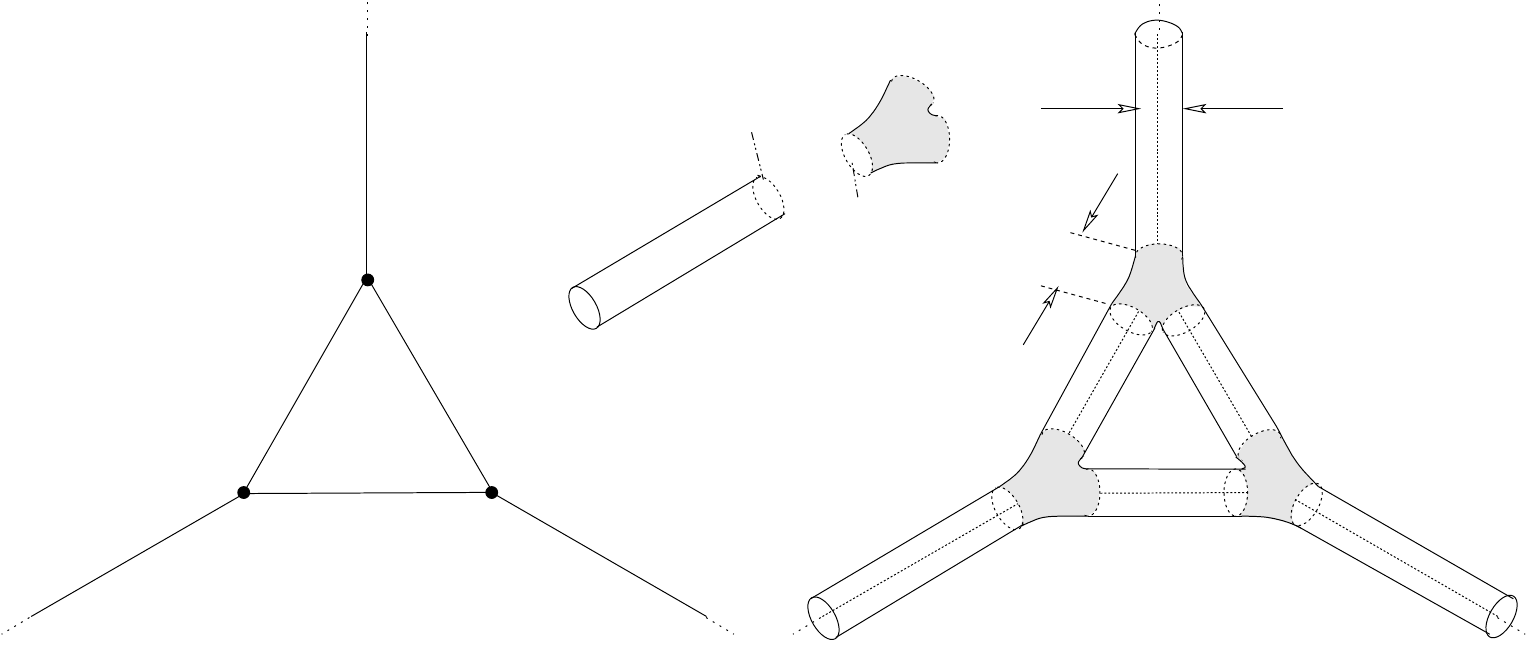}%
  \end{picture}%
  \setlength{\unitlength}{2072sp}
  \begin{picture}(11385,4891)(263,-4360)
    \put(9426,-16){$\sim \eps$}%
    \put(8050,-1426){$\sim \eps$}%
    \put(811,-3481){$e$}%
    \put(1996,-3466){$v$}%
    \put(2566,-3976){$X_0$}%
    \put(6426,-3346){$X_\edeps$}%
    \put(8011,-3700){$X_\vxeps$}%
    \put(8706,-4176){$X_\eps$}%
    \put(4836,-2201){$X_\edeps$}%
    \put(7626,-511){$X_\vxeps$}%
    \put(5436,-156){$Y_\edeps$}%
    \put(6581,-1161){$Y_\edeps$}%
  \end{picture}%
  \caption{A part of a metric graph, the (scaled) building blocks and
    the corresponding (part of a) thin branched manifold (here,
    $Y_e=\Sphere^1$, and $X_\eps$ is the surface of the pipeline
    network). The vertex neighbourhoods $X_\vxeps$ are drawn in gray.}
  \label{fig:example}
\end{figure}
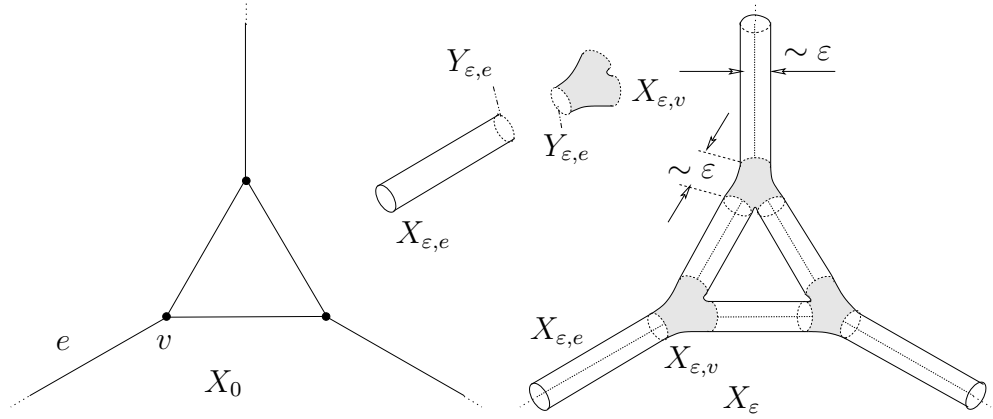

\begin{remark}
  \indent
  \begin{enumerate}[(a)]
  \item Note that we write $X_\vxeps$ etc.\ for a \emph{Riemannian}
    manifold.  More precisely, we should write $(X_v,g_\vxeps)$ for
    this manifold, since the underlying space can be chosen to be
    $\eps$-independent, the $\eps$-dependence only enters via the
    metric.

  \item The manifold $X_\eps$ may have boundary or not. If $X_\eps$
    has boundary, then also (some of) the transversal manifolds $Y_e$
    have boundary.

  \item We define the space $X_\eps$ in an abstract manner, although
    we have concrete examples in mind.  If we consider a graph $M$
    embedded in, say, $\R^d$, and if $\wt X_\eps$ denotes its
    $\eps$-neighbourhood, then we can define a similar decomposition
    as in~\eqref{eq:decomp}, but the building blocks $\wt X_\edeps$
    and $\wt X_\vxeps$ are only \emph{approximatively} isometric with
    $M_e \times \eps Y_e$ and $\eps X_v$ for some fixed Riemannian
    manifolds $(Y_e,h_e)$ and $(X_v,g_v)$.  This may have two reasons:
    \begin{itemize}
    \item We need a little space for the vertex neighbourhoods (of
      order $\eps$), so that we need to replace the interval $M_e$ by
      a slightly smaller one of length $\ell_e - \Err(\eps)$.
    \item The edges may be embedded as non-straight \emph{curves} in
      $\R^2$.  This leads to a slight deviation from the product
      metric.
    \end{itemize}
    All these cases can be treated as a \emph{perturbation} of the
    abstract situation above, see e.g.~\cite[Sec.~5.4 and
    Sec.~6.7]{post:12}).
  \end{enumerate}
\end{remark}

\subsubsection*{Laplacians on thin branched manifolds}

The underlying Hilbert space is $\Lsqr {X_\eps}$ (with the natural
volume measure denoted by $\dd X_\eps$ induced by the $\eps$-dependent
Riemannian metric).  In particular, we have
\begin{multline*}
  \normsqr[\Lsqr {X_\eps}] u
  = \int_{X_\eps} \abssqr{u(x)} \dd X_\eps(x)\\
  = \eps^m \sum_{e \in E} \int_0^{\ell_e} \int_{Y_e} \abssqr{u_e(s_e,y)}
         \dd Y_e(y) \dd s_e
  + \eps^{m+1} \sum_{v \in V} \int_{X_v} \abssqr{u_v(x)} \dd X_v(x)
\end{multline*}
using the decomposition~\eqref{eq:decomp} and suitable
identifications.  In particular, $u_e$ resp.\ $u_v$ denote the
restriction to $X_\edeps$ and $X_\vxeps$, respectively.

As operator on $X_\eps$, we consider the Laplacian $\Delta_\eps \ge 0$
(with Neumann boundary conditions if $\bd X_\eps \ne \emptyset$).
This operator is given as
\begin{equation*}
  (\Delta_\eps u)_v = \frac 1{\eps^2} \Delta_{X_v} u_v
  \quadtext{and}
  (\Delta_\eps u)_e = -u_e'' + (\id \otimes \Delta_{Y_\edeps}) u_e,
\end{equation*}
where $\Delta_{Y_\edeps} \phi=\eps^{-2} \Delta_{Y_e} \phi$ for a
smooth function $\phi$ on $Y_e$, and $\Delta_{Y_e}$ is the (Neumann)
Laplacian on $Y_e$.  Moreover, $(\cdot)'_e$ denotes the derivative
with respect to the longitudinal variable $s \in M_e$.

\subsection{Convergence of the Laplacian on thin branched manifolds}
Let us first define a suitable identification operator
\begin{equation*}
  \map {J_\eps} {\Lsqr {X_0}}{\Lsqr {X_\eps}}.
\end{equation*}

For simplicity, we assume here that $\vol_m (Y_e,h_e)=1$.
As identification operator we choose
\begin{equation*}
  (J_\eps f)_e = f_e \otimes \1_\edeps
  \qquad\text{and}\qquad
  (J_\eps f)_v = 0
\end{equation*}
where $(J_\eps f)_e$ is the contribution on the edge neighbourhoood
$X_\edeps$ and $(J_\eps f)_v$ is the contribution on the vertex
neighbourhood, according to the decomposition~\eqref{eq:decomp}.
Moreover, $\1_\edeps$ is the constant function on $Y_\edeps$ with
value $\eps^{-m/2}$ (the first normalised eigenfunction of
$Y_\edeps$).

\begin{remark}
  The setting $(J_\eps f)_v=0$ seems at first sight a bit rough, but we
  cannot set something like $(J_\eps f)_v=\eps^{-m/2}f(v)$, since on
  $\Lsqr {X_0}$, the value of $f$ at $v$ is not defined.  There is a
  finer version of identification operators on the level of the
  quadratic form domains, again see~\cite[Ch.~4]{post:12} for details.
\end{remark}

Let us now calculate the resolvent difference $R_\eps J_\eps - J_\eps R_0$: For
$g \in \Lsqr{X_0}$ and $w \in \Lsqr{X_\eps}$, we have
\begin{align*}
  \iprod[\Lsqr{X_\eps}] {(R_\eps J_\eps - J_\eps R_0)g} w &=
  \iprod[\Lsqr{X_\eps}] {J_\eps g} {R_\eps w} -
  \iprod[\Lsqr{X_\eps}] {J_\eps R_0g} w\\
  &= \iprod[\Lsqr{X_\eps}] {J_\eps \Delta_0 f} u -
  \iprod[\Lsqr{X_\eps}] {J_\eps f} {\Delta_\eps u},\\
  \intertext{where $u=R_\eps w \in \dom \Delta_\eps$ and $f=R_0 g \in
    \dom \Delta_0$.  Moreover, using the definition of $J_\eps f$, we
    obtain}%
  \dots
  &= \sum_{e \in E} \Bigl( \iprod[\Lsqr{X_\edeps}] {(-f_e''
    \otimes \1_\edeps} {u_e} - \iprod[\Lsqr{X_\edeps}] {f_e \otimes
    \1_\edeps} {-u_e''}\\
  & \hspace{0.3\textwidth}
    - \iprod[\Lsqr{X_\edeps}] {f_e \otimes
    \1_\edeps} {(\id \otimes \Delta_{Y_\edeps}) u_e}
  \Bigr)\\
  &= \sum_{e \in E} \Bigl( \iprod[\Lsqr{X_\edeps}] {(-f_e'' \otimes
    \1_\edeps} {u_e} - \iprod[\Lsqr{X_\edeps}] {f_e \otimes \1_\edeps}
  {-u_e''} \Bigr)\\ %
  \intertext{since we can bring $(\id \otimes \Delta_{Y_\edeps})$ on
    the other side of the inner product (the operator is
    self-adjoint!) and $\Delta_{Y_\edeps} \1_\edeps=0$.  Using $\dd
    X_\edeps= \eps^m \dd Y_e \dd s$ and performing a partial
    integration (Green's first formula), we obtain}%
  \dots
  &= \sum_{e \in E} \eps^{m/2} \Bigl[ \int_{Y_e} (-f_e' \conj u_e +
  f_e \conj u_e') \dd Y_e \Bigr]_{\bd M_e}.\\
  \intertext{Using the conventions $f_e(v)=f_e(0)$ resp.\ $f_e(v)=f_e(\ell_e)$,
  $u_e(v):=u_e(0,\cdot)$ resp. $u_e(v)=u_e(\ell_e,\cdot)$  and $f_e'(v)=
  -f_e'(0)$ resp.\ $f_e'(v)=f_e'(\ell_e)$ if $v=\bd_-e$ resp.\ $v=\bd_+ e$,
  and after reordering, we obtain}%
  \dots
  &= \sum_{e \in E}
  \sum_{v=\bd_\pm e} \eps^{m/2} \int_{Y_e} \bigl(-f_e'(v) \conj u_e(v)
  +
  f_e(v) \conj u_e'(v)\bigr) \dd Y_e\\
  &= \sum_{v \in V} \sum_{e \in E_v} \eps^{m/2} \int_{Y_e}
  \bigl(-\underbrace{f_e'(v) \conj u_e(v)}_{=: I_1} +
        \underbrace{f_e(v) \conj u_e'(v)}_{=: I_2} \bigr) \dd Y_e.
\end{align*}
Consider now
\begin{equation*}
  \avint_v u_v := \frac 1 {\vol X_v} \int_{X_v} u_v \dd X_v
  \qquad\text{and}\qquad
  \avint_e u_e(v) := \frac 1 {\vol Y_e} \int_{Y_e} u_e(v) \dd Y_e,
\end{equation*}
then we can express the sums over $I_1$ as
\begin{align*}
  \sum_{e \in E_v} \eps^{m/2} \int_{Y_e}
  f_e'(v) \conj u_e(v)
  &= \sum_{e \in E_v} \eps^{m/2}
      f_e'(v) \bigl(\avint_e \conj u_e(v)- \avint_v \conj u_v  \bigr)
   +  \Bigl(\sum_{e \in E_v} \eps^{m/2}
      f_e'(v) \Bigr) \avint_v \conj u_v\\
  &= \sum_{e \in E_v} \eps^{m/2}
      f_e'(v) \bigl(\avint_e \conj u_e(v)- \avint_v \conj u_v \bigr).
\end{align*}
The last sum in the first line vanishes since $f \in \dom \Delta_0$
fulfils the so-called Kirchhoff condition
$\sum_{e \in E_v} f_e'(v)=0$.  For the second summand $I_2$, we use
the fact that $f_e(v)=f(v)$ is independent of $e \in E_v$, hence
\begin{align*}
  \sum_{e \in E_v} \eps^{m/2} \int_{Y_e} f_e(v)
                             \conj u_e'(v) \bigr) \dd Y_e
  &= \eps^{m/2} f(v) \int_{\bd X_v} \normder \conj u_v \dd \bd X_v\\
  &= \eps^{m/2} f(v) \int_{X_v} \Delta_{X_v} \conj u_v \dd X_v,
\end{align*}
performing again a partial integration (Green's first formula, writing
$u_v$ as $1 \cdot u_v$).  Summing up the contributions, we have
\begin{align*}
  \iprod[\Lsqr{X_\eps}] {(R_\eps J_\eps - J_\eps R_0)g} w
  &= \sum_{v \in V} \eps^{m/2}
  \bigl(
    - \sum_{e \in E_v} f'_e(v)
        \bigl(\avint_e \conj u_e(v)- \avint_v \conj u_v \bigr)
    + f(v) \int_{X_v} \Delta_{X_v} \conj u_v \dd X_v
  \bigr)\\
  &=: -\iprod[\VXmax]{B_0 g}{A_\eps w}
     + \iprod[\VX] {A_0 g} {B_\eps w},
\end{align*}
where $\VX:=\lsqr{V,\deg}$ (with norm given by
$\normsqr[\lsqr{V,\deg}] \phi := \sum_{v \in V} \abssqr{\phi(v)} \deg
v < \infty$), $\VXmax := \bigoplus_{v \in V} \C^{E_v}$ and
\begin{align*}
  \map {B_0} {&\Lsqr{X_0}} \VXmax,&
  (B_0 g)_v &= \bigl((R_0 g)'_e(v)\bigr)_{e \in E_v},\\
  \map {A_\eps} {&\Lsqr{X_\eps}} \VXmax,&
  (A_\eps w)_v &= \eps^{m/2}
             \bigl( \avint_e (R_\eps w)_e(v) - \avint_v (R_\eps w)_v
                        \bigr)_{e \in E_v}\\
  \map {B_\eps} {&\Lsqr{X_\eps}} \VX,&
  (B_\eps w)(v) &= \frac{\eps^{m/2}}{\deg v}
             \int_{X_v} \Delta_{X_v} (R_\eps w) \dd X_v\\
  \map {A_0} {&\Lsqr{X_0}} \VX,&
  (A_0 g)(v) &= (R_0 g)(v).
\end{align*}
In particular, we have shown
\begin{proposition}
  \label{prp:res.expr}
  We can express the resolvent differences of $\Delta_\eps$ and
  $\Delta_0$, sandwiched with the identification operator $J_\eps$, as
  \begin{equation*}
    \map{R_\eps J_\eps - J_\eps R_0 = -A_\eps^*B_0 + B_\eps^* A_0}
    {\Lsqr{X_0}} {\Lsqr{X_\eps}}
  \end{equation*}
\end{proposition}

We will now show that the $\eps$-dependent operators are actually
small if $\eps \to 0$.  In order to do so, we need two important
estimates:
\subsubsection*{Sobolev trace estimate}
We have
\begin{equation}
  \label{eq:sob.est.mfd}
  \normsqr[\Lsqr{Y_e}] {u(0,\cdot)}
  \le C(\ell_0) \normsqr[\Sob{X_{v,e}}] u
  \bigl(\normsqr[\Lsqr{X_{v,e}}] u
  + \normsqr[\Lsqr{X_{v,e}}] {\nabla u}
  \bigr)
\end{equation}
for all $\map u {X_{v,e}} \C$ smooth enough, where
$X_{v,e}=[0,\ell_0] \times Y_e$ is a collar neighbourhood of the
boundary component of $X_v$ touching the edge neighbourhood $X_e$.
The optimal constant is actually $C(\ell_0)=\coth(\ell_0/2)$.  The
proof of~\eqref{eq:sob.est.mfd} ist just a vector-valued version of
the Sobolev trace estimate
\begin{equation}
  \label{eq:sob.est.int}
  \abssqr{f(0)} \le C(\ell_0)
  \bigl(\normsqr[{\Lsqr{[0,\ell_0]}}] f
  + \normsqr[{\Lsqr{[0,\ell_0]}}] {f'}
  \bigr).
\end{equation}

\subsubsection*{A min-max estimate}
We have
\begin{equation}
  \label{eq:min.max.est}
  \normsqr[\Lsqr{X_v}] {u - \avint u}
  \le \frac 1 {\lambda_2(X_v)} \normsqr[\Lsqr{X_v}] {\nabla u}
\end{equation}
for all $u$ smooth enough, where $\lambda_2(X_v)$ is the first
(non-vanishing) Neumann eigenvalue of $X_v$.  Note that $u - \avint u$
is the projection onto the space orthogonal to the first (constant)
eigenfunction on $X_v$.

As a consequence, we obtain
\begin{lemma}
  \label{lem:av.int}
  We have
  \begin{equation*}
    \eps^m \sum_{e \in E_v} \bigabssqr{\avint_e u_e(v) - \avint_v u}
    \le \eps C(\ell_0)
         \Bigl(\frac 1 {\lambda_2(X_v)} + 1 \Bigr)
               \normsqr[\Lsqr{X_\vxeps}] {\nabla u}.
  \end{equation*}
\end{lemma}
\begin{proof}
  We have (denoting by $\ell_0>0$ a lower bound on the edge lengths)
  \begin{align*}
    \eps^m \sum_{e \in E_v} \bigabssqr{\avint_e u_e(v) - \avint_v u}
    &= \eps^m \sum_{e \in E_v} \bigabssqr{\avint_e (u - \avint_v u)}\\
    &\le \eps^m \sum_{e \in E_v} \int_{Y_e} \bigabssqr{u - \avint_v u}
                         \dd Y_e\\
    &\le \eps^m C(\ell_0) \sum_{e \in E_v}
         \Bigl(\normsqr[\Lsqr{X_{v,e}}] {u - \avint_v u}
               + \normsqr[\Lsqr{X_{v,e}}] {\nabla u}
         \Bigr)\\
    &\le \eps^m C(\ell_0)
         \Bigl(\normsqr[\Lsqr{X_v}] {u - \avint_v u}
               + \normsqr[\Lsqr{X_v}] {\nabla u}
         \Bigr)\\
    &\le \eps C(\ell_0)
         \Bigl(\frac 1 {\lambda_2(X_v)} + 1 \Bigr)
               \normsqr[\Lsqr{X_\vxeps}] {\nabla u}\\
  \end{align*}
  using Cauchy-Schwarz in the first inequality, \eqref{eq:sob.est.mfd}
  and the fact that $\nabla \avint_v u=0$ in the second estimate, the
  fact that $\bigcup_{e \in E_v} X_{v,e} \subset X_v$ in the third
  estimate and~\eqref{eq:min.max.est} and the scaling behaviour
  $\eps^{m-1} \normsqr[\Lsqr{X_v}] {\nabla u}=\normsqr[\Lsqr{X_\vxeps}]
  {\nabla u}$ in the fourth estimate.
\end{proof}

The following result is not hard to see using the Sobolev trace
estimate~\eqref{eq:sob.est.int}:
\begin{proposition}
  \label{prp:ab.0}
  Assume that $0 < \ell_0 \le \ell_e$ for all $e \in E$, then the
  operators $A_0$ and $B_0$ are bounded by a constant depending only
  on $\ell_0$.
\end{proposition}

\begin{proposition}
  \label{prp:ab.eps}
  Assume that
  \begin{equation}
    \label{eq:main-ass}
    0 < \ell_0 \le \ell_e \; \forall e \in E, \qquad
    0 < \lambda_2 \le \lambda_2(X_v)
    \qquad\text{and}\qquad
    \frac {\vol X_v}{\deg v} \le c_{\vol} < \infty  \; \forall v \in V
  \end{equation}
  holds, then $\norm{A_\eps} = \Err(\eps^{1/2})$ and
  $\norm{B_\eps} = \Err(\eps^{3/2})$, and the errors depend only on
  $\ell_0$, $\lambda_0$ and $c_{\vol}$.
\end{proposition}
\begin{proof}
  For $A_\eps$, we have
  \begin{align*}
    \normsqr[\VXmax]{A_\eps w}
    &= \eps^m \sum_{v \in V} \sum_{e \in E_v}
    \bigabssqr{\avint_e u_e(v) - \avint_v u_v}\\
    &\le \eps C(\ell_0)
         \Bigl(\frac 1 {\lambda_2} + 1 \Bigr)
               \sum_{v \in V} \normsqr[\Lsqr{X_\vxeps}] {\nabla u}
    \le \eps C(\ell_0)
         \Bigl(\frac 1 {\lambda_2} + 1 \Bigr)
               \normsqr[\Lsqr{X_\eps}] {\nabla u}
  \end{align*}
  using \Lem{av.int}, where $u=R_\eps w$.  Now, since $u \in \dom
  \Delta_{X_\eps}$, and since $\Delta_{X_\eps}$ is the operator
  associated with the quadratic form, we have
  \begin{equation*}
    \label{eq:spec.calc}
    \normsqr[\Lsqr{X_\eps}]{\nabla u}
    = \iprod[\Lsqr{X_\eps}]{\Delta_{X_\eps} u} u
    = \iprod[\Lsqr{X_\eps}]{\Delta_{X_\eps} (\Delta_{X_\eps}+1)^{-1}w}
      {(\Delta_{X_\eps}+1)^{-1}w}
    \le \normsqr[\Lsqr{X_\eps}] w
  \end{equation*}
  and the inequality is true by the spectral calculus.

  For $B_\eps$, we have
  \begin{align*}
    \normsqr[\VX] {B_\eps g}
    = \eps^m \sum_{v \in V}
       \frac 1 {\deg v} \Bigabssqr{\int_{X_v} \Delta_{X_v} u}
    &\le \eps^m \sum_{v \in V} \frac{\vol X_v}{\deg v}
          \normsqr[\Lsqr {X_v}] {\Delta_{X_v} u}\\
    &= \eps^3  \sum_{v \in V} \frac{\vol X_v}{\deg v}
          \normsqr[\Lsqr {X_\vxeps}] {\Delta_{X_\vxeps} u}\\
    &\le \eps^3  c_{\vol}
     \normsqr[\Lsqr {X_\eps}] {\Delta_{X_\eps} (\Delta_{X_\eps} + 1)^{-1} w)}
    \le \eps^3  c_{\vol}
     \normsqr[\Lsqr {X_\eps}] w
  \end{align*}
  using the scaling behaviour $\Delta_{X_\vxeps}=\eps^{-2}
  \Delta_{X_v}$ and $\normsqr[\Lsqr{X_\vxeps}] w = \eps^{m+1}
  \normsqr[\Lsqr{X_v}] w$, where again $u=R_\eps w$.
\end{proof}

Combining the previous results (\Prpss{res.expr}{ab.0}{ab.eps}),
we have shown the following:
\begin{theorem}
  \label{thm:main1}
  Assume that~\eqref{eq:main-ass} holds, then
  \begin{equation*}
    \norm[\Lsqr{X_0} \to \Lsqr{X_\eps}] {R_\eps J_\eps - J_\eps R_0}
    = \Err(\eps^{1/2}),
  \end{equation*}
  where the error depends only on $\ell_0$, $\lambda_0$ and $c_{\vol}$.
\end{theorem}

\begin{theorem}
  Assume that~\eqref{eq:main-ass} holds, then the (Neumann) Laplacian
  $\Delta_{X_\eps}$ converges to the standard (Kirchhoff) Laplacian
  $\Delta_{X_0}$ in the generalised norm resolvent sense.

  In particular, the results of \Thm{que} apply, i.e., we have
  convergence of the spectrum (discrete or essential) and we can
  approximate $\phi(\Delta_{X_\eps})$ by $J_\eps\phi(\Delta_{X_0}) J_\eps^*$ in
  operator norm up to an error of order $\Err(\eps^{1/2})$.
\end{theorem}
\begin{proof}[Idea of proof]
  We have to show that $J_\eps$ is $\delta_\eps$-quasi unitary.
  It is not hard to see that
  \begin{equation*}
    (J_\eps^* u)_e(s)  =  \eps^{m/2} \int_{Y_e} u_e(s,\cdot)\dd Y_e,
  \end{equation*}
  and that
  \begin{subequations}
    \begin{equation}
      \label{eq:que-1-ex}
      J_\eps^* J_\eps f = f
    \end{equation}
    for all $f \in \Lsqr{X_0}$ (i.e., going from the metric graph to
    the manifold and back, we do not loose information).

    Hence we only have to show that
    \begin{multline}
      \label{eq:que-2-ex}
      \normsqr{u - J_\eps J_\eps^* u}
      = \sum_{v \in V} \normsqr[\Lsqr{X_\vxeps}] {u_v}
      + \sum_{e \in E} \int_{M_e} \normsqr[\Lsqr{Y_\edeps}]
      {u_e(s,\cdot) - \avint_e u_e(s,\cdot)} \dd s\\
      \le \delta_\eps^2 \normsqr{(\Delta_\eps + 1) u}
    \end{multline}
  \end{subequations}
  for some $\delta_\eps \to 0$.  Actually, this can be done using
  similar ideas as before.  For details, we refer again
  to~\cite[Sec.~6.3]{post:12}, and one can show that $\delta_\eps =
  \Err(\eps^{1/2})$ under the additional assumption that $0<\lambda_2'
  \le \lambda_2(Y_e)$ (the first non-zero eigenvalue of $\Delta_{Y_e}$
  on $Y_e$).
\end{proof}

\begin{remark}
  \label{rem:grieser}
  Note that Grieser showed in~\cite{grieser:08} that the $k$-th
  eigenvalue $\lambda_k(\Delta_\eps)$ of the (Neumann) Laplacian
  converges to the $k$-th eigenvalue of the metric graph (Kirchhoff)
  Laplacian $\lambda_k(\Delta_0)$, i.e.,
  \begin{equation*}
    \lambda_k(\Delta_\eps)-\lambda_k(\Delta_0)=\Err(\eps),
  \end{equation*}
  for compact metric graphs and a corresponding family of compact
  Riemannian manifolds using asymptotic expansions.  From our
  analysis, we only obtain the error $\Err(\eps^{1/2})$ as we use less
  elaborated methods.  From Grieser's result, it follows that
  $\dSPEC(R_\eps,R_0)=\Err(\eps)$ where $R_\eps=(\Delta_\eps+1)^{-1}$.
  We conclude from~\eqref{eq:dist.que-uni} that
  $\dQUE(R_\eps,R_0)=\Err(\eps)$.  Our identification operator
  $J_\eps$ only shows the estimate
  $\dQUE(R_\eps,R_0) \le \Err(\eps^{1/2})$.  Knowing already
  Grieser's result, we can directly define a unitary map sending
  eigenfunctions of the metric graph to eigenfunctions of the
  manifold.

  Our identification operator $J_\eps$ just immitates the
  eigenfunctions up to an error.  It would be interesting to see
  whether one can also obtain the optimal error estimate $\Err(\eps)$
  using other identification operators $J_\eps$ respecting in more
  detail the domains and also the local structure of the spaces.
\end{remark}

%
%


\begin{thebibliography}{PS18b}

\bibitem[AP18]{anne-post:pre18}
C.~Ann\'e and O.~Post, \emph{{Wildly perturbed manifolds: norm resolvent and
  spectral convergence}}, arXiv:1802.01124 (2018).

\bibitem[B{\"o}g17]{boegli:17}
S.~B{\"o}gli, \emph{Convergence of sequences of linear operators and their
  spectra}, Integral Equations Operator Theory \textbf{88} (2017), 559--599.

\bibitem[B{\"o}g18]{boegli:18}
\bysame, \emph{Local convergence of spectra and pseudospectra}, J. Spectr.
  Theory \textbf{8} (2018), 1051--1098.

\bibitem[EP05]{exner-post:05}
P.~Exner and O.~Post, \emph{Convergence of spectra of graph-like thin
  manifolds}, Journal of Geometry and Physics \textbf{54} (2005), 77--115.

\bibitem[EP13]{exner-post:13}
\bysame, \emph{A general approximation of quantum graph vertex couplings by
  scaled {S}chr\"odinger operators on thin branched manifolds}, Comm. Math.
  Phys. \textbf{322} (2013), 207--227.

\bibitem[Gri08]{grieser:08}
D.~Grieser, \emph{Spectra of graph neighborhoods and scattering}, Proc. Lond.
  Math. Soc. (3) \textbf{97} (2008), 718--752.

\bibitem[HN99]{herbst-nakamura:99}
I.~Herbst and S.~Nakamura, \emph{{Schr{\"o}dinger operators with strong
  magnetic fields: Quasi-periodicity of spectral orbits and topology}},
  American Mathematical Society. Transl., Ser. 2, Am. Math. Soc.
  \textbf{189(41)} (1999), 105--123.

\bibitem[Kat66]{kato:66}
T.~Kato, \emph{Perturbation theory for linear operators}, Springer-Verlag,
  Berlin, 1966.

\bibitem[KP18]{khrabustovskyi-post:18}
A.~Khrabustovskyi and O.~Post, \emph{Operator estimates for the crushed ice
  problem}, Asymptot. Anal. \textbf{110} (2018), 137--161.

\bibitem[KS03]{kuwae-shioya:03}
K.~Kuwae and T.~Shioya, \emph{Convergence of spectral structures: a functional
  analytic theory and its applications to spectral geometry}, Comm. Anal. Geom.
  \textbf{11} (2003), 599--673.

\bibitem[Lin12]{lin-huaxin:12}
H.~Lin, \emph{Approximate unitary equivalence in simple {$C^\ast$}-algebras of
  tracial rank one}, Trans. Amer. Math. Soc. \textbf{364} (2012), 2021--2086.

\bibitem[Pos06]{post:06}
O.~Post, \emph{Spectral convergence of quasi-one-dimensional spaces}, Ann.
  Henri Poincar\'e \textbf{7} (2006), 933--973.

\bibitem[Pos12]{post:12}
\bysame, \emph{Spectral analysis on graph-like spaces}, Lecture Notes in
  Mathematics, vol. 2039, Springer, Heidelberg, 2012.

\bibitem[PS18a]{post-simmer:18}
O.~Post and J.~Simmer, \emph{Approximation of fractals by discrete graphs: norm
  resolvent and spectral convergence}, Integral Equations Operator Theory
  \textbf{90} (2018), 90:68.

\bibitem[PS18b]{post-simmer:pre18a}
\bysame, \emph{Approximation of fractals by manifolds and other graph-like
  spaces}, arXiv:1802.02998 (2018).

\bibitem[PS19]{post-simmer:pre18d}
\bysame, \emph{Quasi-unitary equivalence and generalised norm resolvent
  convergence}, (preprint) (2019).

\bibitem[RS80]{reed-simon-1}
M.~Reed and B.~Simon, \emph{{Methods of modern mathematical physics I:
  Functional analysis}}, Academic Press, New York, 1980.

\bibitem[Wei00]{weidmann:00}
J.~Weidmann, \emph{Lineare {O}peratoren in {H}ilbertr\"aumen. {T}eil 1},
  Mathematische Leitf\"aden., B. G. Teubner, Stuttgart, 2000, Grundlagen.

\end{thebibliography}
\providecommand{\bysame}{\leavevmode\hbox to3em{\hrulefill}\thinspace}
\providecommand{\MR}{\relax\ifhmode\unskip\space\fi MR }
\providecommand{\MRhref}[2]{%
  \href{http://www.ams.org/mathscinet-getitem?mr=#1}{#2}
}
\providecommand{\href}[2]{#2}



\end{document}